\documentclass[dvips,color,xnumber,francais]{amsart}
\usepackage{cras,math}
\newtheorem*{defnn}{Définition}
\newtheorem*{defnnn}{Définitions}

\begin{document}
\issueinfo{33?}{?}{?}{2000}
\title[Représentations de groupes branchés]{Sous-groupes paraboliques
  et représentations de groupes branchés}
\author{Laurent BARTHOLDI}
\address{Section de Math\'ematiques, Universit\'e de Gen\`eve, CP 240,
  1211 Gen\`eve 24, Suisse}
\email{Laurent.Bartholdi@math.unige.ch}
\author{Rostislav I. GRIGORCHUK}
\address{Steklov Institute of Mathematics, Gubkina 8, Moscow 117966, Russia}
\email{grigorch@mi.ras.ru}
\date{15 mai 2000}
\commby{Mikhael \textsc{Gromov}}
\seriesinfo{Théorie des groupes}{Group Theory}
\thanks{Les auteurs remercient le <<Fonds National Suisse pour la
  Recherche Scientifique>> et l'Université Hébra\"ique de Jérusalem}
\etitle{Parabolic Subgroups and Representations of Branch Groups}
\begin{abstract}
  Soit $G$ un groupe branché (au sens de~\cite{grigorchuk:jibg})
  agissant sur un arbre $\tree$. Un sous-groupe \emdef{parabolique}
  $P$ est le stabilisateur d'un rayon géodésique infini de $\tree$. On
  note $\rho_{G/P}$ la représentation quasi-régulière associée.

  Si $G$ est discret, ces représentations sont irréductibles, mais
  si $G$ est profini, elles se décomposent en une somme directe de
  représentations de dimension finie
  $\rho_{G/P_{n+1}}\ominus\rho_{G/P_n}$, où $P_n$ est le stabilisateur 
  d'un sommet de niveau $n$ de $\tree$.

  Pour quelques exemples concrets, on décompose complètement
  $\rho_{G/P_n}$ en composantes irréductibles.  $(G,P_n)$
  et $(G,P)$ sont des paires de Gelfand, d'où de nouvelles
  occurences d'algèbres de Hecke abéliennes.
\end{abstract}
\begin{eabstract}
  Let $G$ be a branch group (in the sense of~\cite{grigorchuk:jibg})
  acting on a tree $\tree$. A \emdef{parabolic} subgroup
  $P$ is the stabiliser of an infinite geodesic ray in $\tree$. We
  denote by $\rho_{G/P}$ the associated quasi-regular representation.
  
  If $G$ is discrete, these representations are irreducible, but if
  $G$ is profinite, they split as a direct sum of finite-dimensional
  representations $\rho_{G/P_{n+1}}\ominus\rho_{G/P_n}$, where $P_n$
  is the stabiliser of a level-$n$ vertex in $\tree$.

  For a few concrete examples, we completely split
  $\rho_{G/P_n}$ in irreducible components.  $(G,P_n)$
  and $(G,P)$ are Gelfand pairs, whence new occurrences of abelian
  Hecke algebra.
\end{eabstract}
\maketitle

\begin{englishversion}
  We initiate in this paper the study of unitary representations of
  groups of branch and fractal nature. The importance of this class of
  groups was made clear in~\cite{grigorchuk:jibg,bartholdi:phdfr} and
  the numerous papers cited therein.
  
  A \emdef{branch group} is a group $G$ acting ``nicely'' on a
  $d$-regular rooted tree $\tree$, with a finite-index subgroup $K$
  such that $K^d$ embeds in $K$ by acting on the $d$ subtrees below
  the root of $\tree$. We introduce parabolic subgroups $P$, which are
  stabilizers of infinite rays in $\tree$, and establish that they are
  weakly maximal. We then study the corresponding quasi-regular
  representations $\rho_{G/P}$.
  
  If $G$ is a discrete branch group, these representations are
  irreducible. If $G$ is a profinite branch group, $\rho_{G/P}$ is a
  direct sum of the trivial representation and of the
  finite-dimensional representations
  $\rho_{G/P_{n+1}}\ominus\rho_{G/P_n}$, where $P_n$ is the stabilizer
  of a level-$n$ vertex in the tree on which $G$ acts.
  
  We formulate our results by focusing explicitly on a few typical
  examples, one of which, $\overline\Gamma$, is the first occurrence
  of a torsion-free branch group. The facts known of our examples are
  assembled in the following table:
  \def\0{$+$} \def\1{$-$}
  \begin{center}
    \begin{tabular}{r|ccccccc|}
      & $\Gg$ & $\tilde\Gg$ & $\Gamma$ & $K$ & $\overline\Gamma$ &
      $\overline K$ & $\doverline\Gamma$\\
      & \hspace{1cm} & \hspace{1cm} & \hspace{1cm} & \hspace{1cm} &
      \hspace{1cm} & \hspace{1cm} & \hspace{1cm}{ }\\[-12pt]\hline
      Just-infinite &               \0 &   \0 &   \0 & \0 &   \1 & \1 &   \0\\
      Just-nonsolvable &            \0 &   \0 &   \0 & \0 &   \0 & \0 &   \0\\
      Branch &                      \0 &   \0 &   \0 & \0 &   \1 & \1 &   \0\\
      Weak Branch &                 \0 &   \0 &   \0 & \0 &   \0 & \0 &   \0\\
      Fractal &                     \0 &   \0 &   \0 & \1 &   \0 & \0 &   \0\\
      Congruence Property &         \0 &   \0 &   \0 & \0 &   ?  & ?  &   \0\\
      Torsion &                     \0 &   \1 &   \1 & \1 &   \1 & \1 &   \0\\
      (1)\hfill Virtually Torsion-free&\1& \1 &   \0 & \0 &   \0 & \0 &   \1\\
      (2)\hfill Intermediate growth&\0 &   \0 &   \0 & \0 &   \0 & \0 &   \0\\
      (3)\hfill Finitely $L$-Presented&\0& \0 &   \0 & \0 &   \0 & \0 &   \0\\
      (4)\hfill Finite Width&       \0 &   \0 &   \0 & \0 &   ? & ? &     \1\\[4pt]
      \hline
    \end{tabular}
  \end{center}
  
  Property~(1) ranks among the main contributions of this note. (2) is
  studied in~\cite{grigorchuk:growth,bartholdi:phd}.  (3) is studied
  in~\cite{lysionok:pres,bartholdi-g:parabolic,bartholdi:lpres}. (4),
  conjectured for $\Gg$ as early as 1989 by the second author, is
  studied in~\cite{rozhkov:lcs,bartholdi-g:lie,bartholdi:lcs}. A
  thorough treatment of $\Gg$ appears in~\cite{harpe:cgt}.

  For these examples, we completely split the quasi-regular
  representations $\rho_{G/P_n}$ in irreducible components.  $(G,P_n)$
  and $(G,P)$ are Gelfand pairs, giving new instances of abelian Hecke
  algebra.
\end{englishversion}

\section{Introduction}
Nous amorçons dans cet article l'étude des représentations
unitaires de groupes fractals et branché. L'importance de cette
classe de groupes a été mise en évidence
dans~\cite{grigorchuk:jibg,bartholdi:phdfr} et les nombreuses
références qui s'y trouvent.

Nous définissons des sous-groupes paraboliques $P$, démontrons
leur maximalité faible, et étudions les représentations
quasi-régulières $\rho_{G/P}$ associées.

Si $G$ est un groupe branché discret, ces représentations sont
irréductibles. Si $G$ est un groupe branché profini, $\rho_{G/P}$
est la somme directe de la représentation triviale de $G$ et des
représentations de dimension finie
$\rho_{G/P_{n+1}}\ominus\rho_{G/P_n}$, où $P_n$ est le stabilisateur
d'un sommet de niveau $n$ dans l'arbre sur lequel $G$ agit.

Nous avons formulé nos résultats en nous concentrant sur quelques
exemples typiques ; l'un d'eux, $\overline\Gamma$, est le premier
exemple de groupe branché virtuellement sans torsion. Nous donnons
plus de détails dans~\cite{bartholdi-g:parabolic}.

Nous décomposons complètement les représentations $\rho_{G/P_n}$
en composantes irréductibles. Pour ces exemples, $(G,P_n)$ et $(G,P)$
sont des paires de Gelfand, produisant de nouvelles occurrences
d'algèbres de Hecke abéliennes.

\section{Définitions principales}
Les groupes que nous considérons sont tous des sous-groupes du
groupe $\aut(\tree)$ d'automorphismes d'un arbre régulier $\tree$.
Soit $\Sigma=\{1,\dots,d\}$ un alphabet fini. L'ensemble des sommets
de l'arbre $\tree_\Sigma$ est l'ensemble des suites finies sur
$\Sigma$; deux suites sont reliées par une arête si une suite peut
être obtenue à partir de l'autre par ajout à droite d'une lettre
de $\Sigma$. La racine de l'arbre est la suite vide $\emptyset$, et
les descendants du sommet $\sigma$ sont tous les $\sigma s$, pour
$s\in\Sigma$.

Choisissons un ensemble $\Sigma$, et soit $\tree=\tree_\Sigma$. Pour
tout sous-groupe $G<\aut(\tree)$, soit $\stab_G(\sigma)$ le
sous-groupe de $G$ constitué des automorphismes fixant la suite
$\sigma$, et soit $\stab_G(n)$ le sous-groupe de $G$ constitué
des suites fixant toutes les suites de longueur $n$:
\[\stab_G(\sigma)=\{g\in
G|\,g\sigma=\sigma\},\qquad\stab_G(n)=\bigcap_{\sigma\in\Sigma^n}\stab_G(\sigma).\]
Les $\stab_G(n)$ sont des sous-groupes distingués d'indice fini dans
$G$ ; en particulier $\stab_G(1)$ est d'indice au plus $d!$. Soit
$G_n$ le quotient $G/\stab_G(n)$. Si $g\in\aut(\tree)$ fixe la suite
$\sigma$, on note $g_{|\sigma}$ l'élément
correspondant de $\aut(\tree)$ associé à la restriction aux
suites commençant par $\sigma$, ce qui s'écrit en formules $\sigma
g_{|\sigma}(\tau)=g(\sigma\tau)$.  Comme le sous-arbre en-dessous de
n'importe quel sommet est isomorphe à l'arbre initial
$\tree_\Sigma$, on obtient ainsi une application
\begin{equation}
  \psi:\begin{cases}\stab_{\aut(\tree)}(1)\to\aut(\tree)^\Sigma\\h\mapsto (h_{|1},\dots,h_{|d})\end{cases}
\end{equation}
qui est un isomorphisme de groupes.

Un sous-groupe $G<\aut(\tree)$ est \emph{sphériquement transitif} si
l'action de $G$ sur $\Sigma^n$ est transitive pour tout $n\in\N$. On
supposera toujours que cette condition est satisfaite.

$G$ est \emdef{fractal} si pour chaque sommet $\sigma$ de
$\tree_\Sigma$ on a $\stab_G(\sigma)_{|\sigma}\cong G$, où
l'isomorphisme est induit par l'identification de $\tree_\Sigma$ avec
le sous-arbre enraciné en $\sigma$. On a alors pour tout $n$ une
injection $\stab_G(n)<G^{d^n}$.

Si $\sigma$ est une suite et $g\in\aut(\tree)$ est un automorphisme, on
note $g^\sigma$ l'élément de $\aut(\tree)$ agissant comme $g$ sur
les suites commençant par $\sigma$, et trivialement sur les aures:
$g^\sigma(\sigma\tau)=\sigma g(\tau)$, et $g^\sigma(\tau)=\tau$ si
$\tau$ ne commence pas par $\sigma$.

Soit $G<\aut(\tree)$ un groupe agissant fidèlement et
sphériquement transitivement sur un arbre enraciné $\tree_\Sigma$.
Le \emdef{stabilisateur rigide} de $\sigma$ est
$\rist_G(\sigma)=\{g^\sigma|\,g\in G\}\cap G$, et on note
$\rist_G(n)=\prod_{\sigma\in\Sigma^n}\rist_G(\sigma)$.

\begin{defnnn}
  Soit $G<\aut(\tree)$ un groupe sphériquement transitif.
  \begin{enumerate}
  \item $G$ est \emdef{régulièrement branché} s'il possède un
    sous-groupe $K<\stab_G(1)$ d'indice fini tel que
    $K^\Sigma<\psi(K)$.
    \label{defn::rb}
  \item $G$ est \emdef{branché} si $\rist_G(n)$ est d'indice fini dans
    $G$ pour tout $n$. \label{defn::b}
  \item $G$ est \emdef{faiblement branché} si tous ses
    stabilisateurs rigides $\rist_G(\sigma)$ sont infinis.
    \label{defn::wb}
  \item $G$ est \emdef{rugueux} si tous ses stabilisateurs rigides
    sont finis.
  \end{enumerate}
\end{defnnn}
Remarquons que, pour les groupes fractals, \ref{defn::rb} implique
\ref{defn::b} implique \ref{defn::wb} dans la définition ci-dessus, et
qu'un groupe est soit faiblement branché soit rugueux. Remarquons
aussi qu'en principe ces notions dépendent du choix de $\tree$ et de
l'action de $G$.

Toutes ces définitions sont aussi valables dans la catégorie des
groupes profinis, mais dans ce cas il faut considérer $\aut(\tree)$
comme un groupe profini avec sa toplogie naturelle, et $G$ doit être
un sous-groupe fermé.

\section{Exemples principaux}
Comme exemple de groupes rugueux, il y a $\Z$, $D_\infty=\Z/2*\Z/2$ et
le groupe de l'allumeur de réverbères $\Z/2\wr\Z$. On se focalise
ici plutôt sur les groupes branchés et faiblement branchés.

\subsection{Le groupe $\Gg$}\label{subs:defG}
Ce groupe a été défini par le second auteur en
1980~\cite{grigorchuk:burnside}, et agit sur l'arbre $\tree_2$. Soit
$a$ l'automorphisme de $\tree_2$ permutant les deux branches à la
racine. Soit récursivement $b$ l'automorphisme agissant comme $a$
sur la branche gauche et comme $c$ sur la branche droite, $c$
l'automorphisme agissant comme $a$ à gauche et comme $d$ à droite,
et $d$ l'automorphisme agissant trivialement à gauche et comme $b$
à droite. En formules, $\psi(b)=(a,c)$, $\psi(c)=(a,d)$ et
$\psi(d)=(1,b)$.  Soit $\Gg$ le groupe engendré par $\{a,b,c,d\}$.
(Un quelquonque des générateurs $\{b,c,d\}$ peut être omis, car
$\{1,b,c,d\}$ est le groupe de Klein.)

\subsection{Le groupe $\tilde\Gg$}
Cet autre groupe a été défini par les auteurs
dans~\cite{bartholdi-g:lie}, et agit aussi sur $\tree_2$. Avec la même
notation que ci-dessus, on définit $\tilde b,\tilde c,\tilde d$ par
$\psi(\tilde b)=(a,\tilde c)$, $\psi(\tilde c)=(1,\tilde d)$ et
$\psi(\tilde d)=(1,\tilde b)$.  Soit $\tilde\Gg$ le groupe engendré
par $\{a,\tilde b,\tilde c,\tilde d\}$.  Clairement, tous ces
générateurs sont d'ordre $2$ et $\{\tilde b,\tilde c,\tilde d\}$
engendre un groupe abélien élémentaire d'ordre $8$.  Aussi,
$\tilde\Gg$ a pour sous-groupe $\Gg=\langle a,\tilde b\tilde c,\tilde
c\tilde d,\tilde d\tilde b\rangle$. Soit $\tilde{\mathfrak K}$ le
sous-groupe normal $\langle[a,\tilde b],[a,\tilde
d]\rangle^{\tilde\Gg}$ de $\tilde\Gg$.

\subsection{Le groupe $\Gamma$}
Les troix groupes suivants sont des sous-groupes de $\aut(\tree_3)$.
Soit $a$ l'automorphisme de $\tree_3$ permutant cycliquement les trois
branches supérieures. Soit $r$ l'automorphisme de $\tree_3$ défini
récursivement par $\psi(r)=(a,1,r)$.  Soit $\Gamma$ le sous-groupe de
$\aut(\tree_3)$ engendré par $\{a,r\}$. Soit $K$ le sous-groupe normal
$\langle ar,ra\rangle$ de $\Gamma$.

\subsection{Le groupe $\overline\Gamma$}
Soit $s$ l'automorphisme de $\tree_3$ défini récursivement par
$\psi(s)=(a,a,s)$, et soit $\overline\Gamma$ le sous-groupe de
$\aut(\tree_3)$ engendré par $\{a,s\}$. Soit $\overline K$ le
sous-groupe normal $\langle sa^{-1},a^{-1}s\rangle$ de
$\overline\Gamma$.

\subsection{Le groupe $\doverline\Gamma$}
Soit $t$ l'automorphisme de $\tree_3$ défini récursivement par
$\psi(t)=(a,a^{-1},t)$, et soit $\doverline\Gamma$ le sous-groupe de
$\aut(\tree_3)$ engendré par $\{a,t\}$; il a été étudié dans
les années 80 par Narain Gupta and Said
Sidki~\cite{gupta-s:3group,gupta-s:infinitep}.

\section{Propriétés algébriques}
Nous résumons les propriétés principales de nos exemples dans le
tableau ci-dessous. Remarquons que tous ces exemples sont
résiduellement finis. Un point d'interrogation (?) indique que la
propriété n'est pas connue pour ce groupe.

Rappelons qu'un groupe $G$ est \emdef{juste-infini} s'il est infini,
mais que tous les quotients propres sont finis. $G$ est
\emdef{juste-non-résoluble} s'il n'est pas résoluble, mais que
tous ses quotients propres le sont.  $G$ a la \emdef{propriété de
  congruence} si tout sous-groupe de $G$ d'indice fini contient
$\stab_G(n)$ pour un certain $n$.  $G$ a \emdef{croissance
  intermediaire} si pour une quelconque métrique des mots sur $G$ le
volume des boules croît à un taux plus rapide que polynômial
mais plus lent qu'exponentiel.  $G$ est de \emdef{largeur finie} s'il
existe une borne uniforme sur le rang des quotients
$\gamma_n(G)/\gamma_{n+1}(G)$ de sa suite centrale descendante.  $G$
est de \emdef{$L$-présentation finie} s'il peut être présenté
par un ensemble fini $S$ de générateurs et les itérées d'un
ensemble fini de relations par un ensemble fini de substitutions sur
$S$.

\def\0{$+$} \def\1{$-$}
\def\2{\hspace{20pt}}
\begin{center}
\begin{tabular}{r|ccccccc|}
& $\Gg$ & $\tilde\Gg$ & $\Gamma$ & $K$ & $\overline\Gamma$ &
$\overline K$ & $\doverline\Gamma$\\
& \2 & \2 & \2 & \2 & \2 & \2 & \2{ }\\[-12pt]\hline
Juste-infini &                \0 &   \0 &   \0 & \0 &   \1 & \1 &   \0\\
Juste-non-résoluble &       \0 &   \0 &   \0 & \0 &   \0 & \0 &   \0\\
Régulièrement branché sur & {\small $\langle[a,b]\rangle^\Gg$} &
$\tilde{\mathfrak K}$ & $\Gamma'$ & $\Gamma'$ & \1 & \1 &
$\doverline\Gamma'$\\
Faiblement branché &        \0 &   \0 &   \0 & \0 &   \0 & \0 &   \0\\
Fractal &                     \0 &   \0 &   \0 & \1 &   \0 & \0 &   \0\\
Propriété de congruence & \0 &   \0 &   \0 & \0 &   ? & ? &     \0\\
Torsion &                     \0 &   \1 &   \1 & \1 &   \1 & \1 &   \0\\
(1)\hfill\ Sous-groupe sans torsion d'indice&$\infty$& $\infty$ & 3 & 1 & 3 & 1 & $\infty$\\
(2)\hfill\ Croissance intermédiaire&\0 &   \0 &   \0 & \0 &   \0 & \0 &   \0\\
(3)\hfill\ $L$-présentation finie&\0& \0 &   \0 & \0 &   \0 & \0 &   \0\\
(4)\hfill\ Largeur finie&       \0 &   \0 &   \0 & \0 &   ? & ? &     \1\\[4pt] \hline
\end{tabular}
\end{center}

La propriété~(1) fait partie des contributions majeures de cette
note. (2) est étudiée dans~\cite{grigorchuk:growth,bartholdi:phdfr}.
(3) est étudiée
dans~\cite{lysionok:pres,bartholdi-g:parabolic,bartholdi:lpres}. (4),
conjecturée pour $\Gg$ dès 1989 par le second auteur, est
étudiée dans~\cite{rozhkov:lcs,bartholdi-g:lie,bartholdi:lcs}. La
référence~\cite{harpe:cgt} dévoue un chapitre entier à $\Gg$.

\section{Sous-groupes paraboliques}
Soit $\tree=\Sigma^*$ un arbre enraciné. Un \emdef{rayon} $e$ dans
$\tree$ est une géodésique infinie partant de la racine de $\tree$,
ou de manière équivalente un élément de
$\partial\tree=\Sigma^\N$.

Soit $G<\aut(\tree)$ un sous-groupe quelconque et soit $e$ un rayon.
Le \emdef{sous-groupe parabolique} associé est $P_e=\stab_G(e)$.

Les points suivants méritent une attention particulière :
\begin{itemize}
\item $\bigcap_{e\in\partial\tree}P_e=\bigcap_{g\in G}P^g=1$
  (c'est-à-dire que $P$ a un \emdef{noyau trivial}).
\item Soit $e=e_1e_2\dots\in\Sigma^\N$ un rayon infini est soient
  $P_n=\stab_G(e_1\dots e_n)$ les sous-groupes stabilisant un point de
  niveau $n$. Alors les $P_n$ sont d'indice $d^n$ dans $G$ (puisque
  $G$ agit transitivement sur $\Sigma^n$) et satisfont
  $P_e=\bigcap_{n\in\N}P_n$.
\item $P$ est d'indice infini dans $G$, et a la même image que $P_n$
  dans le quotient $G_n$.
\end{itemize}

Soit $G$ un groupe branché, et soit $H$ un sous-groupe quelconque.
$H$ est \emdef{faiblement maximal} si $H$ est d'indice infini dans
$G$, mais si tous les sous-groupes de $G$ contenant strictement $H$
sont d'indice fini dans $G$.  (On remarque que tout sous-groupe infini
de type fini possède des sous-groupes faiblement maximaux, par le
lemme de Zorn).

\begin{thm}
  Soit $P$ un sous-groupe parabolique d'un groupe régulièrement
  branché $G$. Alors $P$ est faiblement maximal.
\end{thm}

Si $G$ est un groupe branché, le sous-groupe parabolique $P$ peut
être décomposé explicitement en une extension scindée
itérée par des groupes finis. Par exemple, pour le groupe
prototypique $\Gg$ du paragraphe~\ref{subs:defG} pour lequel $d=2$, on
pose $e=dd\dots$ et $P=P_e$, obtenant le
\begin{thm}
  $P/P'$ est un $2$-groupe infini élémentaire, engendré par les
  images de $c$, $d$ et des éléments de la forme
  $(1,\dots,1,(ac)^4)$ dans $\rist_\Gg(n)$ pour tout $n\in\N$. On a la
  décomposition suivante:
  \[P = \bigg(B\times \Big(\left(K\times
    \left((K\times\dots)\rtimes\langle
      (ac)^4\rangle\right)\right)\rtimes\langle
  b,(ac)^4\rangle\Big)\bigg)\rtimes\langle c,(ac)^4\rangle,\] où
  chaque facteur $B,K,K,\dots$ de niveau $n$ dans la décomposition
  agit sur le sous-arbre juste sous $e_n$ mais ne contenant pas
  $e_{n+1}$.  $B$ est ici la clôture normale de $b$, et $K$ est la
  clôture normale de $[a,b]$.
\end{thm}

\section{Représentations quasi-régulières}
Si $H$ est un sous-groupe du groupe discret $G$, on note $\rho_{G/H}$
la représentation quasi-régulière de $G$ sur $\ell^2(G/H)$; si
$H=1$, elle est la représentation régulière gauche de $G$.

Le \emdef{commensurateur} du sous-groupe $H$ de $G$ est
\[\comm_G(H) = \{g\in G|\,H\cap H^g\text{ est d'indice fini dans 
  }H\text{ et }H^g\}.\]

De façon équivalente, $\comm_G(H)$ est l'ensemble des $g\in G$
tels que les classes à droite $gH$ et $g^{-1}H$ ont des orbites
finies dans $\{kH\}_{k\in G}$ pour l'action de $H$ sur $G/H$ par
multiplication à gauche.

Un critère de George Mackey affirme que, pour un groupe infini $G$,
la représentation quasi-régulière $\rho_{G/H}$ est
irréductible si et seulement si $\comm_G(H)=H$.

\begin{thm}
  Si $G$ est un groupe discret fractal et faiblement branché, alors
  $\comm_G(P)=P$.
\end{thm}

\begin{cor}
  Il y a une quantité non dénombrable de représentations
  irréductibles non-équivalentes de la forme $\rho_{G/P}$, où
  $G$ est fractal, faiblement branché, et $P$ est un sous-groupe
  parabolique.
\end{cor}

Une situation complètement différente se présente si $G$ est un
groupe profini ; soient en effet $\rho_{G/P_n}$ les représentations
de dimension finie. Elles forment une tour ascendante de
représentations, avec
$\rho_{G/P_{n+1}}=\rho_{G/P_n}\oplus\pi_n^\perp$ pour des
représentations $\pi_n^\perp$. On remarque aussi que le sous-groupe
$P=\bigcap_{n\ge0}P_n$ est fermé.
\begin{thm}
  Soit $G$ un groupe profini branché et soit $P$ un sous-groupe
  parabolique.  Alors la représentation quasi-régulière
  $\rho_{G/P}$ se décompose en une somme de représentations de
  dimension finie :
  \[\rho_{G/P} = \bigoplus_{n\ge0}\pi_n^\perp.\]
\end{thm}

Supposons maintenant que ce groupe profini branché est la complétion
profinie $\widehat G$ d'un groupe discret $G$ satisfaisant la
propriété de congruence. Alors les représentations
irréductibles de $\widehat G$ sont en correspondence bi-univoque
avec les représentations irréductibles des groupes finis
$G_n=G/\stab_G(n)$. Or $\rho_{G_n}$ est une sous-représentation de
$\rho_{G/P_n}\otimes\dots\otimes\rho_{G/P_n}$ (avec $d^n$ facteurs),
puisque $\stab_G(n)=\bigcap_{g\in G/P_n}P_n^g$. On est ainsi amené
à étudier des représentations $\rho_{G/P_n}$ et de leurs
composantes irréductibles. On poursuit cette étude sur nos
exemples de base dans la section suivante.

\section{Algèbres de Hecke et paires de Gelfand}
Considérons maintenant les représentations quasi-régulières
$\rho_{G/P_n}$ de dimension finie, qui se factorisent à travers le
groupe fini $G_n$.

\begin{defnn}
  Soit $G$ un groupe et $H$ un sous-groupe.  L'\emdef{algèbre de
    Hecke} (aussi appelée \emdef{algèbre d'intersection})
  $\Hecke(G,H)$ est l'algèbre $\operatorname{End}_G(\ell^2(G/H))$,
  où $\ell^2(G/H)$ est un $G$-module pour la multiplication à
  gauche.  $\Hecke(G,H)$ peut être pensée comme l'algèbre des
  fonctions $(H,H)$-biinvariantes sur $G$, avec le produit de
  convolution
  \[(f\cdot g)(x) = \sum_{y\in G} f(xy) g(y^{-1}).\]
\end{defnn}

$\Hecke(G,H)$ est engendrée par les fonctions $(H,H)$-biinvariantes
sur $G$, où de façon équivalente par les doubles classes $HgH$.
Le résultat suivant souligne l'importance des algèbres de Hecke
dans l'étude de la décomposition des
représentations~\cite[Section~11D]{curtis-r:methods} : supposons que
$H$ est d'indice fini dans $G$. Alors $\Hecke(G,H)$ est une algèbre
semi-simple.  Il y a une bijection canonique entre les composantes
irréductibles de $\rho_{G/H}$ et les facteurs simples de
$\Hecke(G,H)$, et préservant leurs multiplicités.

Ainsi, si $\Hecke(G,P)$ est abélienne, sa décomposition en
$G$-modules simples a autant de composantes qu'il y a de doubles
classes $PgP$ dans $G$.  Les doubles classes de $P_n$ dans $G$, elles,
sont données par les orbites de $P_n$ sur $G/P_n$, ou, en d'autres
termes, par les orbites de $P_n$ sur $\Sigma^n$. Celles-ci peuvent
être décrites explicitement :

\begin{lem}
  Pour les groupes $\Gg$ et $\tilde\Gg$, les sous-groupes $P_n$ et
  $\tilde P_n$ ont $n+1$ orbites dans $\Sigma^n$; ce sont $\{2^n\}$ et
  les $2^i1\Sigma^{n-1-i}$ pour $0\le i<n$. Les orbites des groupes
  profinis $P$ et $\tilde P$ dans $\Sigma^\N$ sont $\{2^\infty\}$ et
  les $2^i1\Sigma^\N$ pour tout $i\in\N$.
  
  Pour les trois exemples $\Gamma$, $\overline\Gamma$ et
  $\doverline\Gamma$, le sous-groupe $P_n$ a $2n+1$ orbites dans
  $\Sigma^n$; ce sont $\{3^n\}$ et les $3^i1\Sigma^{n-1-i}$ et
  $3^i2\Sigma^{n-1-i}$ pour $0\le i<n$. Les orbites des groupes
  profinis $P$ dans $\Sigma^\N$ sont $\{3^\infty\}$ et les
  $3^i1\Sigma^\N$ et $3^i2\Sigma^\N$ pour tout $i\in\N$.
\end{lem}

\begin{thm}
  $\rho_{\Gg/P_n}$ et $\rho_{\tilde\Gg/\tilde P_n}$ se décomposent
  en somme directe de $n+1$ composantes irréductibles, une de
  degré $2^i$ pour tout $i\in\{1,\dots,n-1\}$ et deux de degré
  $1$.
  
  $\rho_{\Gamma/P}$, $\rho_{\overline\Gamma/P}$ et
  $\rho_{\doverline\Gamma/P}$ se décomposent en somme directe de
  $2n+1$ composantes irréductibles, deux de degré $3^i$ pour tout
  $i\in\{1,\dots,n-1\}$ et trois de degré $1$.
\end{thm}

On a vu que l'algèbre de Hecke $\Hecke(G,P_n)$ est grosso modo de
dimension $n$. Sa structure est encore plus claire si on introduit la
définition suivante : soit $G$ un groupe et $H$ un sous-groupe. La
paire $(G,H)$ est une \emdef{paire de Gelfand} si toutes les
sous-représentations irréductibles de $\rho_{G/H}$ ont
multiplicité $1$.  De façon équivalente, $(G,H)$ est une paire
de Gelfand si $\Hecke(G,H)$ est abélienne.

\begin{thm}
  Pour nos cinq exemples $(G,P)$ et $(G,P_n)$ sont des paires de
  Gelfand pour tout $n\in\N$.
\end{thm}

Ainsi $\Hecke(\Gg,P_n)\cong\Hecke(\tilde\Gg,\tilde P_n)\cong\C^{n+1}$
et $\Hecke(\Gamma,P_n)\cong\Hecke(\overline\Gamma,P_n)\cong
\Hecke(\doverline\Gamma,P_n)\cong\C^{2n+1}$. On a $\Hecke(G,P)\cong\C$
pour tous les groupes discrets considérés, et $\Hecke(\widehat
G,P)\cong\C^\infty$ pour leur complétion profinie.

\bibliography{mrabbrev,people,math,grigorchuk,bartholdi}
\end{document}